\documentclass[a4paper]{article}
\usepackage{amsmath,amsfonts,amssymb}
\usepackage{multicol,multirow,diagbox,bm}
\usepackage{graphicx,subfigure}
\usepackage{lineno,cases}
\allowdisplaybreaks 
\usepackage{pgf,tikz}
\usetikzlibrary{arrows,shapes,automata,positioning}
\usetikzlibrary{fit}
\usepackage{pst-node}
\usepackage{cite}
\usepackage{color,xcolor}
\renewcommand{\bar}{\overline}
\renewcommand{\top}{\intercal}

\newcommand{\e}{\varepsilon}
\newcommand{\R}{\mathbb{R}}
\newcommand{\w}{\omega}

\usepackage{bm} 
\newcommand{\N}{\mathcal{N}}

\newtheorem{ass}{Assumption}

\newtheorem{thm}{Theorem}

\newtheorem{rem}{Remark}

\newcommand{\pb}{\noindent\textbf{Proof.} }
\newcommand{\pe}{\hfill\rule{4pt}{8pt}}

\def\rm{\mathrm}
\usepackage{cite}
\begin{document}
	
\title{Neural Network-based Constrained Optimal Coordination for Heterogeneous Uncertain Nonlinear Multi-agent Systems
	\thanks{This work was supported by National Natural Science Foundation of China under Grants 61973043 and 61773373.}
}
\author{Yutao Tang, Ding Wang
	\thanks{Y. Tang is with the School of Artificial Intelligence, Beijing University of Posts and Telecommunications, Beijing, China.  D. Wang is  with the Faculty of Information Technology, Beijing University of Technology, Beijing, China. (E-mails: yttang@bupt.edu.cn, dingwang@bjut.edu.cn)}
}	
\date{ }
\maketitle

{\noindent\bf Abstract}: In this paper, we investigate a constrained optimal coordination problem for a class of heterogeneous  nonlinear multi-agent systems described by high-order dynamics subject to both unknown nonlinearities and external disturbances.  Each agent has a private objective function and a steady-state constraint about its output.  We develop a composite distributed controller for each agent by a combination of internal model and neural network. All agent outputs are proven to reach the constrained minimal point of the aggregate objective function with bounded residual errors irrespective of the unknown nonlinearities and external disturbances.  Two examples are finally given to demonstrate the effectiveness of the algorithm.


\maketitle

\maketitle

\section{Introduction}

Multi-agent coordination has been a hot topic over the last decades and has many practical applications in multi-robot control, smart grid, and sensor networks \cite{derenick2007convex,stegink2017unifying,zhang2015twc}. As one of the most interesting problems, distributed consensus optimization attracts more and more attention due to the fast development of machine learning and big data technologies. Various effective algorithms have been proposed to achieve such an optimal coordination in different situations \cite{boyd2011distributed,yang2017distributed,nedic2020distributed}. 

Recently, many efforts have been made to incorporate high-order agent dynamics into the distributed optimization design. This is mainly due to the fact that distributed optimization tasks may be implemented or depend on physical plants of high-order dynamics in practice, e.g., source seeking in multi-robot systems \cite{zhang2011extremum} and attitude formation control of rigid bodies \cite{song2017relative}.  
For example, an optimal coordination problem for double integrators was considered in  \cite{zhang2017distributed} with an integral control idea and further extended for Euler-Lagrange agents.  Distributed optimization with bounded controls was also explored for both single and double integrators in  \cite{xie2017global}. For general linear systems, an embedded technique was developed in  \cite{tang2019cyb} to simplify the whole design by converting  the original optimal coordination problem into several subproblems and solving them almost independently. At the same time, some interesting  attempts have also been made for special classes of nonlinear multi-agent systems. For example,  \cite{wang2016distributed} was focused on a class of nonlinear agents in output feedback form with unity relative degree and solved its optimal coordination problem by improving the integral rule in \cite{zhang2017distributed}. The embedded design idea was also further exploited for nonlinear agents in different forms \cite{tang2018distributed, tang2020optimal, liu2022distributed} .  

So far, there are few optimal coordination results considering optimization constraints on the final states except for single-integrator multi-agent systems as mentioned above. Compared with the unconstrained case, the set constraint will pose some specific challenge. In fact, additional mechanisms are usually needed to ensure the satisfaction of constraints on decision variables \cite{jokic2009constrained,glattfelder2012control,tee2009barrier, garone2017reference}. Thus, the design of effective algorithms and the associated convergence analysis are  more involved.  When facing  multi-agent systems of nontrivial dynamics, the problem will inevitably be much more challenging than the optimal coordination results derived either for single integrators or without such constraints. 

Based on these observations, this paper focuses on the optimal coordination problem for a typical class of heterogeneous nonlinear multi-agent systems with set constraints.  Moreover, we assume the high-order agents are subject to both unknown nonlinearities and external disturbances.  To overcome the difficulties brought from the nonlinearity, uncertainty, and constrained optimization requirement, we view the formulated problem as an asymptotic regulation problem where the reference point is determined by the constrained optimization problem and develop a novel neural network-based distributed control to solve the optimal coordination problem. We also provide rigorous theoretical analysis to ensure the global stability of resultant closed-loop systems.   To our knowledge, this might be the first attempt to solve such kind of optimal coordination problems by neural network-based controls in this setting.

The main contributions of this paper is twofold. On the one hand, we present and solve a constrained optimal coordination problem for high-order nonlinear agents. Compared with existing optimal coordination  results \cite{xie2017global,zhang2017distributed, qiu2019distributed,tang2019cyb}, this paper extends them to the case with both set constraints on the global optimization requirement and more general nonlinear agent dynamics. On the other hand, a novel neural network-based controller combined with internal model designs is developed to achieve the optimal coordination goal in a distributed manner. Thanks to the approximation ability of neural networks, the composite design method allows us to handle a large class of nonlinear high-order agents subject to both unknown dynamics and external disturbance generated by certain autonomous linear exosystem.  In fact, by removing the restrictive linearly parameterized condition on unknown nonlinearities, this work explicitly generalizes existing results for linear or nonlinear multi-agent systems \cite{xie2017global,tang2019cyb,qiu2019distributed,tang2018distributed}. 

The remainder of this paper is organized as follows: Problem formulation is presented in Section  \ref{sec:formulation}. Then the main result is provided in Section \ref{sec:main} with detailed designs. Following that, two numerical examples are given to illustrate the efficiency of our algorithm  in Section \ref{sec:simu}. Finally, conclusions are given in Section \ref{sec:con}. 


{\sl Notation}: Let $\R^n$ be the $n$-dimensional Euclidean space and $\R^{n\times m}$ be the set of all $n\times m$ matrices. ${\it  1}_n$ (or ${\it  0}_n$) denotes an $n$-dimensional all-one (or all-zero) column vector and ${\bm 1}_{n\times m}$ (or ${\bm 0}_{n\times m}$) all-one (or all-zero) matrix. 
$\mbox{diag}\{b_1,\,{\dots},\,b_n\}$ denotes an $n\times n$ diagonal matrix with diagonal elements $b_i,\,(i=1,\,{\dots},\,n)$. $\mbox{col}(a_1,\,{\dots},\,a_n) = {[a_1^\top,\,{\dots},\,a_n^\top]}^\top$ for column vectors $a_i\; (i=1,\,{\dots},\,n)$.  For a vector $x$ (or matrix $A$) , $||x||$ ($||A||$) denotes its Euclidean (or spectral) norm. For a square matrix $A$, $\mbox{Tr}(A)$ denotes the trace of $A$ and $||A||_{\rm F}=\mbox{Tr}(A^\top A)$ denotes its Frobenius norm. A continuous function $\alpha\colon[0,\, +\infty)\to [0,\, +\infty)$ belongs to class $\mathcal{K}_\infty$ if it is strictly increasing and satisfies $\alpha(0)=0$ and $\lim_{s\to \infty}\alpha(s)=\infty$.

\section{Problem Formulation}\label{sec:formulation}

In this paper, we consider a collection of $N$ heterogeneous nonlinear systems modeled by:
\begin{align}\label{sys:agent}
{x}^{(n_i)}_{i}&= g_i([x]_i,\, \mu)+b_i u_i+d_i(t),\quad i\in \mathcal{N}=\{1,\,\dots,\,N\}
\end{align}
where $[x]_i\triangleq \mbox{col}(x_{i},\,\dots,\,x^{(n_i-1)}_{i}) \in \R^{q n_i}$ is the state variable of agent $i$ with integer $n_i\geq 2$,  $x_i$ is the output, $u_i \in \mathbb{R}^q$ is the control input, and $\mu\in \R^{n_\mu}$ is an uncertain parameter vector. The high-frequency gain matrix $b_i$ is invertible.  Without loss of generality, we let $b_i=I_q$ and assume the vector-valued function $g_{i}\colon \R^{q n_i}\times \R^{n_\mu}\to \R^q$  to be smooth  but unknown to us. 

The signal $d_i(t)\in \R^q$ represents the external disturbance of agent $i$  which can be modeled by 
\begin{align}\label{sys:disturbance}
\begin{split}
	d_i(t)&=D_i(\mu)\w_i,\quad \dot{\w}_i=S_i\w_i,\quad  \w_i(0)=\w_{i0}\in \R^{m_i}
\end{split}
\end{align}
with $\w_i\in \R^{m_i}$, $S_i\in \R^{m_i\times m_i}$ and  $D_i \in \R^{q\times m_i}$. Moreover, we assume that the matrix  $S_i$ has no eigenvalue with negative real part. 
Note that system \eqref{sys:disturbance} can model many typical disturbances, including a combination of step signals of arbitrary magnitudes, ramp signals of arbitrary slopes, and sinusoidal signals of arbitrary amplitudes and initial phases \cite{huang2004nonlinear}.

As stated in existing publications \cite{kia2015distributed, zhang2017distributed, xie2017global}, we endow this multi-agent system with the following distributed optimization problem
\begin{align}\label{opt:main}
\begin{split}
	&\mbox{minimize} \quad f(y)=\sum\nolimits_{i=1}^{N} f_i(y)\\
	&\mbox{subject to}~~{ y \in \Omega_0\triangleq \cap_{i=1}^N \Omega_i}\subset \R^q
\end{split}
\end{align}
where $f_i\colon \R^q \to \R$ is differentiable.  Furthermore, we assume that each agent only know a part of this optimization problem in the sense that agent $i$ only knows $f_i$ and $\Omega_i$.

To ensure the well-posedness of this optimization problem, we make the following assumption \cite{nesterov2018lectures}.

\begin{ass}\label{ass:convex-strong}
For $i \in \mathcal{N}$,  the set $\Omega_i$ is closed and convex with $\Omega_0$ nonempty;  the function $f_i$ is $\underline{l}_i$-strongly convex and its gradient $\nabla f_i$  is $\overline l_i$-Lipschitz over an open set containing $\Omega_i$ for constants $\underline{l}_i,\, \overline l_i>0$.
\end{ass}

Under this assumption, there exists a unique finite solution to problem \eqref{opt:main} according to Theorem 2.2.10 in \cite{nesterov2018lectures}. Denote it as $y^*=\mbox{argmin}_{y\in \Omega_0} f(y)$.   We aim to regulate the multi-agent system \eqref{sys:agent} such that the agent outputs reach this global minimizer in spite of the uncertainties and disturbances. However, no agent can compute the exact $y^*$ and reach it as expected by itself due the lack of global information of $f$ and $\Omega_0$. In fact, we can introduce some local decision variables for the agents as $y_1,\,\dots,\,y_N$ and denote $\tilde f({\bm y})=\sum_{i=1}^N f_i(y_i)$ with ${\bm y}\triangleq \mbox{col}(y_1,\,\dots,\,y_N)$. Then, the problem \eqref{opt:main} is equivalent to minimize $\tilde f({\bm y})=\sum\nolimits_{i=1}^N f_i(y_i)$ subject to a local set constraint $y_i\in \Omega_i$ and a global consensus constraint $y_1=\dots=y_N$. Since the consensus constraint can only be satisfied via a cooperation, we are more interested in distributed designs where the agents can communicate with some others. 

To this end, we use a directed graph $\mathcal{G}=(\mathcal{N},\, \mathcal{E}, \,\mathcal{A})$ to describe the information sharing relationships among those agents with a node set $\mathcal{N}$, an edge set $\mathcal{E}\subset \mathcal{N}\times \mathcal{N}$, and a weight matrix $\mathcal{A}\in \R^{N\times N}$ \cite{godsil2001algebraic}. If agent $i$ can get the information of agent $j$, then there is an edge $(j,\,i)$ in $\mathcal{E}$, i.e., $a_{ij}>0$.  Here is an assumption to guarantee that any agent's information can reach another. 
\begin{ass}\label{ass:graph}
Graph $\mathcal{G}$ is undirected and connected.
\end{ass}

Regarding multi-agent system \eqref{sys:agent}, function $f_i$, set $\Omega_i$, and graph $\mathcal{G}$,  the constrained optimal coordination problem for agent \eqref{sys:agent} is formulated  to  find a feedback control $u_i$ for agent $i$ by using its own and exchanged information with the neighbors such that all trajectories of agents are well-defined over the time interval $[0,\,+\infty)$ and the resultant outputs satisfy $\lim_{t\to +\infty}\|x_i(t)-y^*\|=0$ for each $i\in \N$.

This optimal coordination problem naturally ensures an output consensus of the multi-agent system \eqref{sys:agent}. Compared with existing output consensus coordination results \cite{ren2008distributed, xi2012output}, the formulation further requires the consensus point to be the optimal solution $y^*$ specified by minimizing a global cost function across the whole network, which is more challenging. It is remarkable that the optimality issue of multi-agent coordination has also been studied from the viewpoint of optimal control \cite{thunberg2016optimal, gao2018leader, wang2019self}. Different from these important results, we emphasize more on the optimal steady-state performance and require the agent outputs reaching a consensus and minimizing some global static optimization problem.

The formulated optimal coordination problem has been partially investigated for second-order agents  \cite{zhang2017distributed, xie2017global, qiu2019distributed}. Here, we further consider  heterogeneous set constraints and higher-order agent dynamics possessing unknown nonlinearities and external disturbances, which inevitability bring more difficulties in achieving such an optimal coordination than these existing results.  

\section{Main Result}\label{sec:main}
To solve this problem, we first consider an auxiliary  problem and convert the original problem into some robust stabilization problem. Then, we develop the final optimal coordination controller by an internal model + neural network-based design.

\subsection{Problem conversion}

We start from the same optimal coordination problem for a group of virtual integrator agents as follows:
\begin{align}\label{sys:agent:virtual}
\dot{r}_i=u_{i}^0
\end{align}
with state $r_i\in \R^q$ and  input $u_{i}^0\in \R^q$. Assign these agents  with the same  cost functions $f_1,\,\dots,\,f_N$ and graph $\mathcal{G}$ as above.  To solve the optimal coordination problem for agent \eqref{sys:agent:virtual} is to develop proper input $u_i^0$ for agent $i$ such that $\lim_{t\to+\infty}r_i(t)=y^*$.  Note that this auxiliary problem has been well-studied in literature and many distributed algorithms can be utilized to solve it, e.g., the ones in  \cite{liu2015second} and  \cite{zeng2017distributed}.

Suppose this auxiliary optimal coordination problem has been done. With these estimates of the global optimal solution $y^*$ given by \eqref{sys:agent:virtual}, we only need to consider a robust tracking problem for agent \eqref{sys:agent} with reference $r_i(t)$ to solve our formulated problem. For better analysis, we denote $x_{i1}=x_i-r_i$, $x_{ij}\triangleq x_{i}^{(j-1)}$ for $2\leq j\leq n_i$. Choose constants $k_{ij}$ for $1\leq j\leq n_i-1$ such that the polynomial $p_i(\lambda)=\sum_{j=1}^{n_i-1}k_{ij} \lambda^{j-1}+\lambda^{n_i-1}$ is Hurwitz. Letting $z_i=\mbox{col}(x_{i1},\,\dots,\, x_{i n_i-1})$ and $\zeta_i=\sum_{j=1}^{n_i-1}k_{ij} x_{ij}+x_{in_i}$ gives an error system as follows:
\begin{align}\label{sys:rd=1}
\begin{split}
	\dot{z}_i&=A_iz_i+ B_i\zeta_i-E_i\dot{r}_i\\
	\dot{\zeta}_i&=\bar g_i(z_i,\,\zeta_i,\, r_i,\, \mu)+u_i+d_i(t)-k_{i1}\dot{r}_i
\end{split}
\end{align}
where $\bar g_i(z_i,\,\zeta_i,\,r_i,\, \mu)=g_{i}([x]_{i},\,\mu)+k_{in_i-1}\zeta_i-k_{in_i-1}k_{i1}x_{i1}+\sum_{j=2}^{n_i-1} (k_{ij-1}-k_{in_i-1}k_{ij})x_{ij}$ and
\begin{align*}
A_i=\left[\begin{array}{c|c}
	{\bm 0}_{n_i-2}&I_{n_i-2}\\\hline
	-k_{i1}&-k_{i2},\dots,-k_{in_i-1}
\end{array}\right]\otimes I_q,\quad B_i=\begin{bmatrix}
	{\bm 0}_{n_i-2}\\
	1
\end{bmatrix}\otimes I_q,\quad E_i=\begin{bmatrix}
	1\\
	{\bm 0}_{n_i-2}
\end{bmatrix}\otimes I_q
\end{align*}


From the above form, we have converted the formulated optimal coordination problem into some robust stabilization problem by viewing $\dot{r}_i$ as a vanishing perturbation.  Compared with similar results for linear systems \cite{xie2017global, tang2019cyb, qiu2019distributed}, our problem involves extra nonlinearities from the set constraints and nonlinear agent dynamics. Moreover, the nonlinearities in this multi-agent system can not be perfectly linearly parameterized as that in  \cite{tang2018distributed}. Consequently, (adaptive) feedback linearization method is not applicable to the associated tracking problem because of the unknown nonlinearity and external disturbances. Then, we have to seek new rules to solve our problem.

Inspired by existing designs \cite{hou2009decentralized,tang2015distributed,tang2020optimal}, we split the whole control effort into two parts as follows:
\begin{align}
u_i=u_{id}+u_{ir}
\end{align}
where $u_{id}$ is designed to compensate the external disturbance and $u_{ir}$ is  to handle the unknown nonlinearity and drive agent $i$ to track its reference $r_i$.

It is well-known that internal model-based control is effective to reject modeled disturbances \cite{huang2004nonlinear}. Here, we construct $u_{id}$ following the same technical line. Let $P_i(s)=s^{n_{p_i}}+\hat{p}_{i1}s^{n_{p_i}-1}+\cdots+\hat{p}_{in_{p_i}-1}s+\hat{p}_{in_{p_i}}$ be the minimal polynomial of  matrix $S_i$ and denote $\tau_i=\text{col}(\tau_{i1},\ldots,\tau_{in_{p_i}})$ with $\tau_{ij}=\frac{{\rm d}^{j-1}d_i(t)}{{\rm d}t^{j-1}}\in \R^q$.  Take two matrices as follows:
\begin{align*}
\Phi_i=\left[\begin{array}{c|c}
	0&I_{n_{p_i}-1}\\\hline
	-\hat{p}_{in_{p_i}}&-\hat{p}_{i{n_{p_i}-1}}\,\cdots\,-\hat{p}_{i1}
\end{array}\right]\otimes I_q,\quad 
\Psi_i=\begin{bmatrix}
	1\\
	{\bm 0}_{n_{pi}-1}
\end{bmatrix}^\top\otimes I_q
\end{align*}
By a direct calculation, we obtain
\begin{equation}\label{dyn:ssg}
\dot{\tau_i}=\Phi_i\tau_i,\quad d_i=\Psi_i \tau_i
\end{equation}
System \eqref{dyn:ssg} is called a steady-state generator \cite{huang2004nonlinear}. Since the pair $(\Psi_i, \Phi_i)$ is observable, there exists a constant matrix  $G_i$ such that $F_i\triangleq \Phi_i+G_i\Psi_i$ is Hurwitz.  To reject the disturbance $d_i$, we propose an internal model-based compensator
\begin{equation}\label{ctr:id}
u_{id}=-\Psi_i\eta_i,\quad \dot{\eta}_i=F_i\eta_i+G_iu_i.
\end{equation}

Next, we are going to propose applicable $u_{ir}$ to complete the whole design.  Since the nonlinear function $g_i([x]_i,\,\mu)$ is unknown to us, the term $\bar g_i(z_i,\,\zeta_i,\,r_i,\,\mu)$ can not be directly used for feedback.  To tackle this issue, an intuitive idea is to estimate this term in some way and develop an estimation-based control law. As neural networks have been proven to be an effective tool to approximate unknown nonlinear functions \cite{lewis1998neural, farrell2006adaptive, ge2002stable, wen2016neuro},  we present a neural network-based rule combined with the above internal model-based compensator to solve our problem in next subsection.

\subsection{Solvability analysis} \label{sec:solvability}
In this subsection, we present the whole design of our optimal coordination rule and provide theoretical stability analysis of the closed-loop system.

For our optimal coordination problem, we aim at global stability performance. Note that neural network-based controls usually ensure control performance in the sense of semiglobal stability of the closed-loop systems \cite{farrell2006adaptive, ge2002stable}. To overcome this shortcoming, we try to utilize the neural networks to estimate the expected feedforwarding control efforts as that in  \cite{chen2012globally}. For this purpose, we let ${\bm u}_i({y}^*,\,\mu)=g_i(\mbox{col}({y}^*,\,{\bm 0}_{q(n_i-1)}),\,\,\mu)$. This is indeed the feedforwarding effort for us to regulate $x_i$ to the optimal point according to Theorem 3.8 in  \cite{huang2004nonlinear}. If the trajectory of $r_i(t)$ indeed converges to the optimal solution $y^*$, it should be uniformly bounded. Thus, we try to reproduce ${\bm u}_i(r_i,\,\mu)$ and develop neural network-based approximation rules for the sequel design. 

To this end, we let $\hat g_i(z_i,\,\zeta_i,\, r_i,\, \mu)\triangleq \bar g_i(z_i,\,\zeta_i,\, r_i,\, \mu)-{\bm u}_i(r_i,\,\mu)$. It is verified that $\hat g_i({\bm 0},\,{\bm 0},\, r_i,\, \mu)={\bm 0}$ for any $r_i\in \R^q$ and $\mu\in \R^{n_\w}$. Motivated by existing neural network-based designs \cite{hou2009decentralized,chen2012globally,ma2015neural}, we use a radial basis function (RBF) network as a function approximator and rewrite ${\bm u}_i$ as follows:
\begin{align*}
{\bm u}_i(r_i,\,\mu)=W_i^\top {\bm \sigma}_i(r_i)+{\bm{\epsilon}}_i(r_i)
\end{align*}
where ${\bm \sigma}_i(r_i)=\mbox{col}({\sigma}_{i1}(r_i),\,\dots,\,{\sigma}_{in_w}(r_i))$ is the activation function vector with ${\sigma}_{ij}(s)=e^{-(s-\mu^{\rm c}_{ij})^2/\kappa_i^2}$ for $j=1,\,\dots,\,n_w$ and $W_i\in\R^{n_w\times q}$ is the weight matrix. Here, $\mu^{\rm c}_{ij}$ is the center of the receptive field, $\kappa_i$ is the width of the Gaussian function, and ${\bm{\epsilon}}_i(r_i)$ is the residual error. By  the universal approximation theorem \cite{farrell2006adaptive}, for any given $\e>0$, there exists an ideal constant weight $W_i^*\in \R^{n_w\times q}$ with a large enough integer $n_w>0$ such that $||{\bm{\epsilon}}_i(r_i)||<\e$ over any compact set. 

Since the ideal weight $W_i^*$ can not be known a prior, we develop the following adaptive neural network-based rule to tackle this issue:
\begin{align}\label{ctr:ir}
\begin{split}
	u_{ir}&=-W_i^\top {\bm \sigma}_i(r_i)-\theta_i\rho_i(\zeta_i)\zeta_i\\ 
	\dot{W}_i&=- \ell (W_i-W^0_{i})+ {\bm \sigma}_i(r_i)\zeta_i^\top   \\
	\dot{\theta}_i&=- \ell(\theta_i-\theta^0_{i})+\rho_i(\zeta_i)||\zeta_i||^2
\end{split}
\end{align}
where function  $\rho_i>0$ is to be specified later. Here $W_i$, $\theta_i$ are dynamic gains, $\ell>0$ is a fixed chosen constant to ensure the boundedness of $W_i$ and $\theta_i$, and the term $-\theta_i \rho_i(\zeta_i)\zeta_i$ is designed to dominate the unknown nonlinearity in \eqref{sys:rd=1}. Similar adaptive controllers have been used in literature \cite{ge2002stable,ma2015neural}. The constants $W_i^0$ and $\theta_i^0$ are chosen parameters based on the (possible) prior information of this multi-agent system, especially the nonlinearities and initial conditions of the whole system. Without further requirements, we can just set the default values as $W_i^0={\bm 0}$ and $\theta_i^0=0$. 

%
%

As for the auxiliary constrained optimal coordination problem for agent \eqref{sys:agent:virtual}, we directly borrow the cooperative laws developed in \cite{liu2015second}. 
Combining \eqref{ctr:id} and \eqref{ctr:ir}, we propose the full optimal coordination controller for agent $i$ as follows:
\begin{align}\label{ctr:full}
\begin{split}
	u_i&=-W_i^\top {\bm \sigma}(r_i)-\theta_i\rho_i(\zeta_i)\zeta_i-\Psi_i\eta_i\\
	\dot{\eta}_i&=F_i\eta_i+G_iu_i\\
	\dot{W}_i&=- \ell (W_i-W^0_{i})+ {\bm \sigma}_i(r_i)\zeta_i^\top \\
	\dot{\theta}_i&=- \ell(\theta_i-\theta^0_{i})+\rho_i(\zeta_i)||\zeta_i||^2\\
	\dot{r}_i&=-2r_i+2P_{\Omega_i} (r_i- \nabla f_i(r_i)-r_i^{\bm o}-v_i^{\bm o})\\
	\dot{v}_i&= r_i
\end{split}
\end{align}
where  $r_i^{\bm o}\triangleq \sum\nolimits_{j=1}^Na_{ij}(r_i-r_j)$, $v_i^{\bm o}\triangleq \sum\nolimits_{j=1}^N a_{ij}(v_i-v_j)$, and $P_{\Omega_i}$ is the projector operator from $\R^q$ to $\Omega_i$. Here, the $\eta_i$ subsystem is the internal model to reject the modeled disturbance \ref{sys:disturbance}, $W_i$, $\theta_i$ are adaptive parameters in the neural networks to approximate the expected feedforwarding control efforts, and $r_i$, $v_i$ are utilized to generate the global optimal solution $y^*$. Clearly, this controller is distributed as agent $i$ only uses its own and neighboring information.

Putting  nonlinear agent \eqref{sys:agent} and distributed controller \eqref{ctr:full} together, we obtain the associated closed-loop system as follows:
\begin{align}
\begin{split}
	\dot{z}_i&=A_iz_i+ B_i\zeta_i-E_i\dot{r}_i\\
	\dot{\zeta}_i&=\hat  g_i(z_i,\,\zeta_i,\, r_i,\, \mu)+{\bm u}_i(r_i,\,\mu)-W_i^\top {\bm \sigma}(r_i)-\theta_i\rho_i(\zeta_i)\zeta_i-\Psi_i\eta_i+d_i-k_{i1}\dot{r}_i \\
	\dot{\eta}_i&=F_i\eta_i+G_iu_i\\
	\dot{W}_i&=- \ell(W_i-W^0_{i})+{\bm \sigma}_i(r_i)\zeta_i^\top\\
	\dot{\theta}_i&=-  \ell (\theta_i-\theta^0_{i})+\rho_i(\zeta_i)||\zeta_i||^2\\
	\dot{r}_i&=-2r_i+2P_{\Omega_i} (r_i- \nabla f_i(r_i)-r_i^{\bm o}-v_i^{\bm o})\\
	\dot{v}_i&= r_i
\end{split}
\end{align}

It is ready to provide our main theorem of this paper.
\begin{thm}\label{thm:main}
Consider the multi-agent system \eqref{sys:agent} with graph $\mathcal{G}$ and function $f_i$ and suppose Assumptions \ref{ass:convex-strong}--\ref{ass:graph} holds. Using controller \eqref{ctr:full} to solve the constrained optimal coordination problem, one has the following results:
\begin{itemize}
	\item[1) ] The trajectories of  $x_i(t),\,\dots,\,x_i^{(n_i-1)}(t)$ are  bounded for all $t\geq 0$, $i=1,\,\dots,\,N$;
	\item[2) ] The coordination errors $||x_i(t)-y^*||$ are uniformly ultimately bounded, i.e., this multi-agent system achieves an approximate optimal coordination with residual errors.
\end{itemize}
\end{thm}
\pb  Let $\bar\eta_i=\eta_i-\tau_i-G_i\zeta_i$, $\bar W_i=W_i-W_i^*$ and $\bar \theta_i=\theta_i-\theta_i^*$ with $\theta_i^*>0$ to be specified later. Then, we have 
\begin{align}\label{sys:error-translated}
\begin{split}
	\dot{z}_i&=A_iz_i+ B_i\zeta_i-E_i\dot{r}_i\\
	\dot{\bar \eta}_i&=F_i\bar \eta_i+F_iG_i\zeta_i-G_i\hat g_i(z_i,\,\zeta_i,\, r_i,\, \mu)+k_{i1}\dot{r}_i\\ 
	\dot{\zeta}_i&=\tilde g_i(z_i,\,\bar \eta_i,\, \zeta_i,\, r_i,\, \mu)-\bar W_i^\top {\bm \sigma}(r_i)-\theta_i^*\rho_i(\zeta_i)\zeta_i -\bar \theta_i\rho_i(\zeta_i)\zeta_i+{\bm{\epsilon}}_i(r_i)-k_{i1}\dot{r}_i \\
	\dot{\bar W}_i&=- \ell (\bar W_i+W_i^*-W^0_{i}) +{\bm \sigma}_i(r_i)\zeta_i^\top\\
	\dot{\bar \theta}_i&=-  \ell (\bar \theta_i+\theta_i^*-\theta^0_{i})+\rho_i(\zeta_i)||\zeta_i||^2\\
	\dot{\bm r}&=- 2{\bm r}+2P_\Omega ({\bm r}- \nabla \tilde f({\bm r})-(L\otimes I_q){\bm r}-(L\otimes I_q){\bm v})\\
	\dot{\bm v}&= {\bm r}
\end{split}
\end{align}
where ${\bm r}=\mbox{col}(r_1,\,\dots,\,r_N)$, ${\bm v}=\mbox{col}(v_1,\,\dots,\,v_N)$, $P_{\Omega}$ is the projector operator determined by $P_{\Omega_i}$, and $\tilde g_i(z_i,\,\bar \eta_i,\,\zeta_i,\,r_i,\, \mu)=\hat g_i(z_i,\,\zeta_i,\, r_i,\, \mu)-\Psi_i\bar \eta_i+\Psi_iG_i\zeta_i$.  It can be easily verified that $\bar g_i({\bm 0},\,{\bm 0},\, r_i,\, \mu)={\bm 0}$ and  $\tilde g_i({\bm 0},\,{\bm 0},\, {\bm 0},\,r_i,\, \mu)=0$ for any $r_i\in \R^q$ and $\mu\in \R^{n_{\mu}}$. 
The proof can be split into two steps as follows:

{\it   Step 1}: We consider the stability of the first two subsystems $\mbox{col}(z_i,\bar \eta_i)$. As the matrices $A_i$ and $F_i$ are Hurwitz, there exist unique positive definite matrices $P_{iz}$ and $P_{i\eta}$ such that $P_{iz}A_i^\top+A_i^\top P_{iz}=-2I_{q(n_i-1)}$ and $P_{i\eta}F_i^\top+F_i^\top P_{i\eta}=-2I_{qn_{p_i}}$. Letting $V_{iz}=z_i^\top P_{iz}z_i$ and $V_{i\bar \eta}=\bar \eta_i^\top P_{i\bar \eta} \bar \eta_i$ gives 
\begin{align*}
\dot{V}_{iz}&=2z_i^\top P_{iz}(A_iz_i+ B_i\zeta_i-E_i\dot{r}_i)\\
&=2z_i^\top P_{iz} A_i z_i+2z_i^\top P_{iz} B_i \zeta_i- 2z_i^\top P_{iz} E_i\dot{r}_i\\
&\leq -||z_i||^2+2||P_{iz} B_i||^2 ||\zeta_i||^2+ 2||P_{iz} E_i||^2||\dot{r}_i||^2
\end{align*}
and 
\begin{align*}
\dot{V}_{i\bar \eta}&=2\bar \eta_i^\top P_{i\eta}[F_i\bar \eta_i+F_iG_i\zeta_i-G_i\hat g_i(z_i,\,\zeta_i,\, r_i,\, \mu)+k_{i1}\dot{r}_i]\\ 
&=2\bar \eta_i^\top P_{i\eta}F_i\bar \eta_i+2\bar \eta_i^\top P_{i\eta} F_iG_i\zeta_i-2\bar \eta_i^\top P_{i\eta}G_i\hat g_i(z_i,\,\zeta_i,\, r_i,\, \mu)+2\bar \eta_i^\top P_{i\eta}k_{i1}\dot{r}_i\\
&\leq -||\bar \eta_i||^2+3 ||P_{i\eta} F_iG_i||^2 ||\zeta_i||^2+ 3||P_{i\eta}k_{i1}||^2||\dot{r}_i||^2 \\&\quad  +3 ||P_{i\eta}G_i||^2 ||\hat g_i(z_i,\,\zeta_i,\, r_i,\, \mu)||^2
\end{align*}

Note that $\hat g_i({\bm 0},\,{\bm 0},\,r_i,\, \mu)=0$ and $\tilde g_i({\bm 0},\,{\bm 0},\, {\bm 0},\,r_i,\, \mu)=0$ for any $r_i$ and $\mu$.  By Lemma 7.8 in  \cite{huang2004nonlinear}, there exist some known smooth functions $\hat \phi_{i1},\hat \phi_{i2}$, $\tilde \phi_{i1},\,\tilde \phi_{i2}>1$ and unknown constants $\hat c_{ig}$, $\tilde c_{ig}>1$ such that 
\begin{align}\label{eq:bound-nonlinear}
\begin{split}
	||\hat g_i(z_i,\,\zeta_i,\, r_i,\, \mu)||^2 &\leq \hat c_{ig}[\hat \phi_{i1}(z_i)||\zeta_i||^2+\hat \phi_{i2}(\zeta_i)||\zeta_i||^2]\\
	||\tilde  g_i(z_i,\,\bar \eta_i,\,\zeta_i,\,r_i,\,\mu)||^2 &\leq \tilde c_{ig}[\tilde \phi_{i1}(\tilde z_i)||\hat z_i||^2+\tilde \phi_{i2}(\zeta_i)||\zeta_i||^2]
\end{split}
\end{align}
where $\tilde z_i\triangleq \mbox{col}(z_i,\,\bar \eta_i)$ for short. 

We apply  Theorem 1 in  \cite{sontag1995changing} to $z_i$-subsystem and obtain that, for any given smooth $\Delta_{iz}(z_i)>0$, there exists a differentiable function $V_{iz}^1(z_i)$ satisfying that
\begin{align*}
&\underline \alpha_{iz}(||z_i||)\leq  V_{iz}^1(z_i)\leq\bar \alpha_{iz}(||z_i||)\\
&\dot{V}_{iz}^1\leq -\Delta_{iz}(z_i) ||z_i||^2+\sigma_{i\zeta}^1\gamma_{i\zeta}^1(\zeta_i)||\zeta_i||^2+\sigma_{ir}^1\gamma_{ir}^1(\dot{r}_i)||\dot{r}_i||^2
\end{align*}
for some known smooth functions $\underline \alpha_{iz},\, \bar \alpha_{iz}\in \cal{K}_\infty$, $\gamma_{i\zeta}^1,\,\gamma_{ir}^1\geq 1$ and unknown constants $\sigma_{i\zeta}^1,\,\sigma_{ir}^1\geq 1$.

Let $V_{i\tilde z}=\ell_{iz} V_{iz}^1(z_i)+V_{i\bar \eta}(\bar \eta_i)$ with a constant $\ell_{iz}>0$ to be specified later. It is positive definite and radially unbounded.  Its time derivative along the trajectory of system \eqref{sys:error-translated} satisfies
\begin{align*}
\dot{V}_{i\tilde z}&\leq -\ell_{iz}\Delta_{iz}(z_i) ||z_i||^2+\ell_{iz}\sigma_{i\zeta}^1\gamma_{i\zeta}^1(||\zeta_i||)||\zeta_i||^2+\ell_{iz}\sigma_{ir}^1\gamma_{ir}(||\dot{r}_i||)||\dot{r}_i||^2 -||\bar \eta_i||^2+3 ||P_{i\eta} F_iG_i||^2 ||\zeta_i||^2\\
&\quad +3 ||P_{i\eta}G_i||^2 ||\hat g_i(z_i,\,\zeta_i,\, r_i,\, \mu)||^2 + 3||P_{i\eta}k_{i1}||^2||\dot{r}_i||^2\\
&\leq -[\ell_{iz}\Delta_{iz}(z_i)-3 \hat c_{ig} ||P_{i\eta}G_i||^2 \hat \phi_{i1}(||z_i||)]||z_i||^2-||\bar \eta_i||^2 +[\ell_{iz}\sigma_{i\zeta}^1\gamma_{i\zeta}^1(\zeta_i)+ 3 ||P_{i\eta} F_iG_i||^2\\
&\quad  + 3 \hat c_{ig} ||P_{i\eta}G_i||^2 \hat \phi_{i2}(\zeta_i)]||\zeta_i||^2 +[\ell_{iz}\sigma_{ir}^1\gamma_{ir}^1(\dot{r}_i)+3||P_{i\eta}k_{i1}||^2]||\dot{r}_{i}||^2
\end{align*}
Let $\Delta_{iz},\,\tilde \gamma_{i\zeta},\,\tilde \gamma_{ir}$ be smooth functions satisfying 
\begin{align*}
\Delta_{iz}(z_i)&\geq  2 \max\{\hat \phi_{i1}(||z_i||),\, 1 \},\quad \tilde \gamma_{i\zeta}(\zeta_i)\geq 3\max\{ \gamma_{i\zeta}^1(\zeta_i)\hat \phi_{i2}(\zeta_i) ,\,1 \},\quad  \tilde \gamma_{ir}(\dot{r}_i)\geq 2 \max\{\gamma_{ir}^1(\dot{r}_i),\,1\}
\end{align*}
and $\ell_{iz},\, \tilde \sigma_{i\zeta},\,\tilde \sigma_{ir}$ be positive constants such that
\begin{align*}
\ell_{iz}&\geq \max\{3 \hat c_{ig} ||P_{i\eta}G_i||^2,\,1\},\quad \tilde \sigma_{i\zeta}\geq \max\{\ell_{iz}\sigma_{i\zeta}^1,\, 3 ||P_{i\eta} F_iG_i||^2\},\quad \tilde \sigma_{ir}\geq\max\{\ell_{iz}\sigma_{ir}^1,\, 3||P_{i\eta}k_{i1}||^2\}
\end{align*}
It follows then
\begin{align*}
\dot{V}_{i\tilde z}&\leq -||\tilde z_i||^2+\tilde \sigma_{i\zeta}\tilde \gamma_{i\zeta}(\zeta_i)||\zeta_i||^2+\tilde \sigma_{ir} \tilde \gamma_{ir}(\dot{r}_i)||\dot{r}_{i}||^2
\end{align*}

{\it   Step 2}: We consider the stability of the $\mbox{col}(\tilde z_i,\, \zeta_i)$-subsystem. Using the changing supply functions technique to this subsystem, one has that, for any given smooth $\Delta_{i\tilde z}(\tilde z_i)>0$, there exists a continuously differentiable function $V_{i\tilde z}^1(\tilde z_i)$ satisfying that
\begin{align*}
&\underline \alpha_{i\tilde z}(||\tilde z_i||)\leq  V_{i\tilde z}^1(\tilde z_i)\leq\bar \alpha_{i\tilde z}(||\tilde z_i||)\\
&\dot{V}_{i\tilde z}^1\leq -\Delta_{i\tilde z}(\tilde z_i) ||\tilde z_i||^2+\tilde \sigma_{i\zeta}^1\tilde \gamma_{i\zeta}^1(\zeta_i)||\zeta_i||^2+\tilde \sigma_{ir}^1\tilde \gamma_{ir}^1(\dot{r}_i)||\dot{r}_i||^2
\end{align*}
for some known smooth functions $\underline \alpha_{i\tilde z},\, \bar \alpha_{i\tilde z}\in \mathcal{K}_\infty$, $\tilde \gamma_{i\zeta}^1,\,\tilde \gamma_{ir}^1\geq 1$ and unknown constants $\tilde \sigma_{i\zeta}^1,\,\tilde \sigma_{ir}^1\geq 1$. 

Let $V_i(\tilde z_i,\,\zeta_i,\bar W_i,\, \zeta_i)=V_{i\tilde z}^1(\tilde z_i)+||\zeta_i||^2+\mbox{Tr}(\bar W_i^\top \bar W_i)+\bar \theta_i^2$. It is positive definite and radially unbounded. Taking its derivative along the trajectory of \eqref{sys:error-translated} gives
\begin{align*}
\dot{V}_i&\leq -\Delta_{i\tilde z}(\tilde z_i) ||\tilde z_i||^2+\tilde \sigma_{i\zeta}^1\tilde \gamma_{i\zeta}^1(\zeta_i)||\zeta_i||^2+\tilde \sigma_{ir}^1\tilde \gamma_{ir}^1(\dot{r}_i)||\dot{r}_i||^2+2\zeta_i^\top[ \tilde g_i(z_i,\,\bar \eta_i,\, \zeta_i,\, r_i)-\bar W_i^\top {\bm \sigma}_i(r_i)-\theta_i^*\rho_i(\zeta_i)\zeta_i-\bar \theta_i\rho_i(\zeta_i)\zeta_i]\\
&\quad +2\zeta_i^\top[{\bm{\epsilon}}_i(r_i)-k_{i1}\dot{r}_i]+2\mbox{Tr}(\bar W_i^\top[- \ell\bar W_i+\ell(W_i^*-W^0_{i})+{\bm \sigma}_i(r_i)\zeta_i^\top])+2\bar \theta_i[-\ell (\bar\theta_i+\theta_i^*-\theta^0_{i})+\rho_i(\zeta_i)||\zeta_i||^2]\\
&\leq -\Delta_{i\tilde z}(\tilde z_i) ||\tilde z_i||^2-2\theta_i^*\rho_i(\zeta_i)||\zeta_i||^2-2\ell \mbox{Tr}(\bar W_i^\top \bar W_i)-2\ell ||\bar \theta_i||^2 +\tilde \sigma_{i\zeta}^1\tilde \gamma_{i\zeta}^1(\zeta_i)||\zeta_i||^2+2\zeta_i^\top \tilde g_i(z_i,\,\bar \eta_i,\, \zeta_i,\, r_i)\\
&\quad -2\zeta_i^\top k_{i1}\dot{r}_i-2\zeta_i^\top {\bm{\epsilon}}_i(r_i) -2\ell\mbox{Tr}(\bar W_i^\top(W_i^*-W^0_{i})) - 2\ell \bar \theta_i^\top (\theta_i^*-\theta_i^0)+\tilde \sigma_{ir}^1\tilde \gamma_{ir}^1(\dot{r}_i)||\dot{r}_i||^2
\end{align*}
where we use the identity $\mbox{Tr}(ab^\top)=b^\top a$ for any two column vectors $a,\, b\in \R^n$.  

Combining this inequality  with \eqref{eq:bound-nonlinear}, we further use Young's inequality and obtain that 
\begin{align*}
\dot{V}_i&\leq -\Delta_{i\tilde z}(\tilde z_i) ||\tilde z_i||^2-2\theta_i^*\rho_i(\zeta_i)||\zeta_i||^2-2\ell \mbox{Tr}(\bar W_i^\top \bar W_i)-2\ell ||\bar \theta_i||^2+\tilde \sigma_{i\zeta}^1\tilde \gamma_{i\zeta}^1(\zeta_i)||\zeta_i||^2+\tilde c_{ig}||\zeta_i||^2\\
&\quad +[\tilde \phi_{i1}(\tilde z_i)||\hat z_i||^2+\tilde \phi_{i2}(\zeta_i)||\zeta_i||^2)]+||\zeta_i||^2+ k_{i1}^2||\dot{r}_i||^2+||\zeta_i||^2+ ||{\bm{\epsilon}}_i(r_i)||^2\\
&\quad - \ell \mbox{Tr}(\bar W_i^\top \bar W_i)- \ell\mbox{Tr}((W_i^*-W^0_{i})^\top(W_i^*-W^0_{i}))+\ell ||\bar \theta_i||^2 +\ell || \theta_i^*-\theta_i^0||^2\\
&\leq -[\Delta_{i\tilde z}(\tilde z_i)-\tilde \phi_{i1}(\tilde z_i)] ||\tilde z_i||^2- \ell \mbox{Tr}(\bar W_i^\top \bar W_i)-[2\theta_i^*\rho_i(\zeta_i)-\tilde \sigma_{i\zeta}^1\gamma_{i\zeta}^1(\zeta_i)-\tilde \phi_{i2}(\zeta_i)-\tilde c_{ig}-2]||\zeta_i||^2\\
&-\ell ||\bar \theta_i||^2 +c_{i\dot{r}}||\dot{r}_i||^2+||{\bm{\epsilon}}_i(r_i)||^2+ \ell||W_i^*-W^0_{i}||_{\rm F}^2+  \ell||\theta_i^*-\theta_i^0||^2
\end{align*}
with $c_{i\dot{r}}\triangleq \sup_{0\leq t\leq+\infty}|\tilde \sigma_{ir}^1 \tilde \gamma_{ir}^1(\dot{r}_i(t))|+k_{i1}^2$. Note that this term $c_{i\dot{r}}$ is well-defined due to the boundedness of $r_i$ and $\dot{r}_i$. 

Choosing $\Delta_{i\tilde z},\, \rho_i$ be smooth functions such  that
\begin{align*}
\Delta_{i\tilde z}(\tilde z_i)\geq 2\max\{\tilde \phi_{i1}(\tilde z_i),\,1\},\quad \rho_i(\zeta_i)\geq \max\{\gamma_{i\zeta}^1(\zeta_i),\, \tilde \phi_{i2}(\zeta_i),\,1 \}
\end{align*}
and $\theta_i^*$ be a constant such that $\theta_i^*\geq \max\{\tilde \sigma_{i\zeta}^1,\,\tilde c_{ig},\,1\}$,
we have 
\begin{align}\label{eq:final}
\dot{V}_i&\leq -||\tilde z_i||^2-||\zeta_i||^2-\ell  \mbox{Tr}(\bar W_i^\top \bar W_i)-\ell ||\bar \theta_i||^2+\Xi_i
\end{align}
where $\Xi_i\triangleq c_{i\dot{r}}||\dot{r}_i||^2+\epsilon^2+ \ell||W_i^*-W^0_{i}||_{\rm F}^2+ \ell ||\theta_i^*-\theta_i^0||^2$.

From the inequality \eqref{eq:final}, we can obtain that the first five subsystems in \eqref{sys:error-translated} is input-to-state stable with $\Xi_i$ as its input by Theorem 4.19  in  \cite{khalil2002nonlinear}. Since $\Xi_{i}$ is upper bounded according to Theorem 2 in   \cite{liu2015second}, we conclude the uniformly ultimate boundedness of trajectories of $\bar \eta_i$ and $\xi_i$ according to Definition 4.7 in  \cite{khalil2002nonlinear}. From the definitions of $\bar \eta_i$, $\xi_i$, and $y^*$, one can obtain the boundedness of state $[x]_i$ and coordination error $x_i(t)-y^*$. This completes the proof.
\pe

\begin{rem}\label{rem:formulation}
This optimal coordination problem has been partially discussed in literature \cite{xie2017global, zhang2017distributed, qiu2019distributed, tang2019cyb} for linear agents. By contrast, the agents here are subject to heterogeneous set constraints and of uncertain nonlinear dynamics. Moreover, 
the developed neural network-based control can also facilitate us to successfully remove the restrictive linearly parameterized condition on nonlinearities required in existing results \cite{tang2018distributed}. 
\end{rem}

\begin{rem}\label{rem:error}
From the expression of $\Xi_{i}$ in inequality \eqref{eq:final},  both $\Xi_{i}$ and the residual error $||x_i(t)-y^*||$ can be made smaller than any given positive constant by selecting a small enough $l$ and increasing the number $n_\w$ of neurons in the neural network. In this sense, this constrained optimal coordination problem for nonlinear multi-agent system \eqref{sys:agent} is solved by our controller \eqref{ctr:full} in a globally practical sense. 
\end{rem}

\section{Numerical Example}\label{sec:simu}

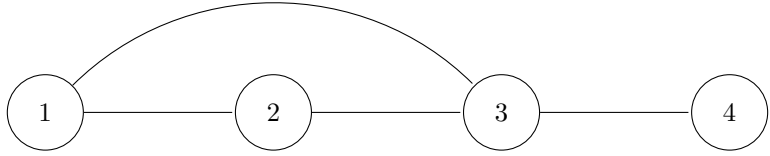
\begin{figure}
\centering
\begin{tikzpicture}[shorten >=1pt, node distance=3 cm, >=stealth',
	every state/.style ={circle, minimum width=1 cm, minimum height=1 cm}]
	\node[align=center,state](node1) {1};
	\node[align=center,state](node2)[right of=node1]{2};
	\node[align=center,state](node3)[right of=node2]{3};
	\node[align=center,state](node4)[right of=node3]{4};
	\path[-]  (node1) edge (node2)
	(node2) edge (node3)
	(node3) edge (node4)
	(node1) edge [bend left=45]  (node3)
	;
\end{tikzpicture}
\caption{Communication graph $\mathcal G$ in our examples.}	\label{fig:graph} 
\end{figure}

In this section, we present two numerical examples to illustrate the effectiveness of our designs.

{\it  Example 1}.  Consider a group of  single-link manipulators with flexible joints \cite{khalil2002nonlinear} modeled by:
\begin{align}\label{sys:manipulator-exam}
J_{i1} \ddot{q}_{i1}+M_igL_i \sin q_{i1}+ k_i(q_{i1}-q_{i2})&=0 \\
J_{i2} \ddot{q}_{i2}-k_i(q_{i1}-q_{i2})&=u_i+d_i,\,i=1,\,\dots,\,4\nonumber
\end{align}
where $q_{i1}, q_{i2}$ are the angular positions, $J_{i1},\, J_{i2}$ ar the moments of inertia, $M_i$ is  the total mass, $L_i$ is a distance, $k_i$ is a spring constant, $u_i$ is the torque input, and $d_i$ is the actuated disturbance of manipulator $i$. Suppose the information sharing graph is depicted in Fig.~\ref{fig:graph} with unity edge weights. It apparently satisfies Assumption \ref{ass:graph}.  Letting $x_i=q_{i1}$, we can rewrite system \eqref{sys:manipulator-exam} into the form \eqref{sys:agent} with $n_i=4$, $b_i=\frac{k_i}{J_{i1}J_{i2}}$ and $g_i([x]_i,\, L_i)=-x_i^{(2)}(\frac{M_igL_i}{J_{i1}}\cos x_i + \frac{k_i}{J_{i1}}+ \frac{k_i}{J_{i2}})+ \frac{M_igL_i}{J_{i1}}(\dot{x}_i^2-\frac{k_i}{J_{i2}})\sin x_i$. 

\begin{figure} 
\centering
\includegraphics[width=0.84\textwidth]{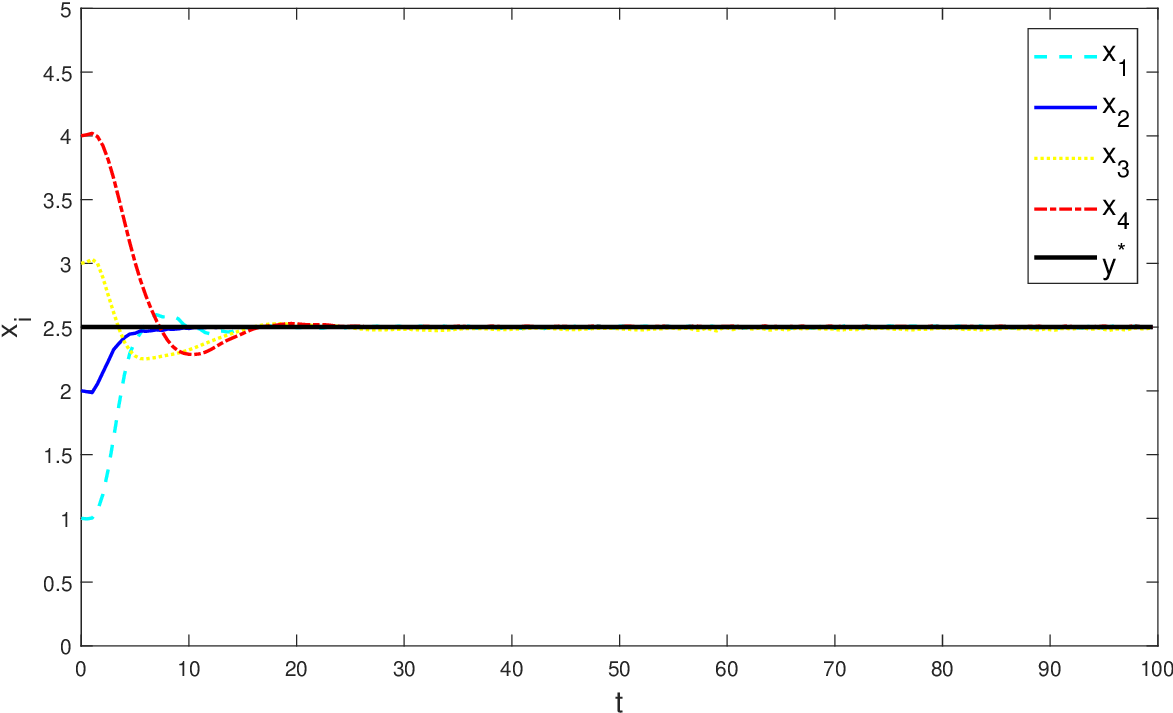}
\caption{Profiles of agent outputs in {\it   Example 1}.}	\label{fig:simu1}
\end{figure}

We want to steer these manipulators to rendezvous at a common position that minimizes the aggregate distance from their starting position to this final position to save resources.  For this purpose, we take the cost functions as $f_i(y)= \frac{1}{2}||y-q_{i1}(0)||^2$ and $f(y)= \frac{1}{2}\sum_{i=1}^4 ||y-q_{i1}(0)||^2$ ($i=1,\,\dots,\, 4$).  To make this problem more interesting, we assume that $L_i=(1+\mu_{i4})L_{i0}$ with nominal length $L_{i0}$ and the external disturbances are described by 
\begin{align*}
D_1&=1+\mu_{15},\quad S_1=0,\quad  D_2=1+\mu_{25},\quad S_2=1\\
D_3&=[1+\mu_{35}~~0],\quad S_3=\begin{bmatrix}
	0&1\\ -1&0
\end{bmatrix},\quad D_4=[1~~1+\mu_{45}~~0],\quad S_4=\mbox{diag}\{1,\,\begin{bmatrix}
	0&2\\-2&0
\end{bmatrix}\}
\end{align*} 
with unknown parameter $\mu_{ij}$.  Here $d_2(t)$ and $d_4(t)$ might tend to infinity depending the initial condition.

Note that feedback linearization rule fails to solve our problem due to the unknown parameters. Nevertheless, we can verify all assumptions in this paper and thus develop a neural network-based control \eqref{ctr:full} for this multi-agent system to solve this problem according to Theorem \ref{thm:main}. To reject those external disturbances for agents, we choose 
\begin{align*}
&F_1=-1,\quad  G_1=-1,\quad F_2=\begin{bmatrix}
	-4&1\\-4&0
\end{bmatrix},\quad G_2=\begin{bmatrix}
	-4\\-4
\end{bmatrix}\\
&F_3=\begin{bmatrix}
	-2&1\\-1&0
\end{bmatrix},\quad G_3=\begin{bmatrix}
	-2\\0
\end{bmatrix},\quad F_4=\begin{bmatrix}
	-4&2\\-2&0
\end{bmatrix},\quad G_4=\begin{bmatrix}
	-4\\0
\end{bmatrix}
\end{align*} for the internal model \eqref{ctr:id}.  In the simulations, we set $J_{i1}=J_{i2}=1$, $L_{i0}=1$, $M_i=1$, $k_i=1$ and assume that the uncertain parameter $\mu_{ij}$ is randomly chosen between $-0.5$ and $0.5$.  To approximate the unknown feedforwarding input,  we construct the RBF neural network with the parameters $n_\w=21$, $\mu_{ij}^{\mbox{c}}=0.5*(j-11)$, and $\kappa_i=1.62$. The nonlinear control gain function is chosen as $\rho_i(s)=s^4+1$ with parameters $k_{i1}=1$, $l=0.01$ for $1\leq i\leq 4$. All initial conditions are randomly chosen. The simulation result is shown in Fig.~\ref{fig:simu1}, where $x_i(t)$ is found to quickly converge to the neighborhood of the optimal point  $y^{\star}=\frac{1}{4}\sum_{i=1}^4 q_{i1}(0)$ with small residual errors.

\begin{figure} 
\centering
\includegraphics[width=0.84\textwidth]{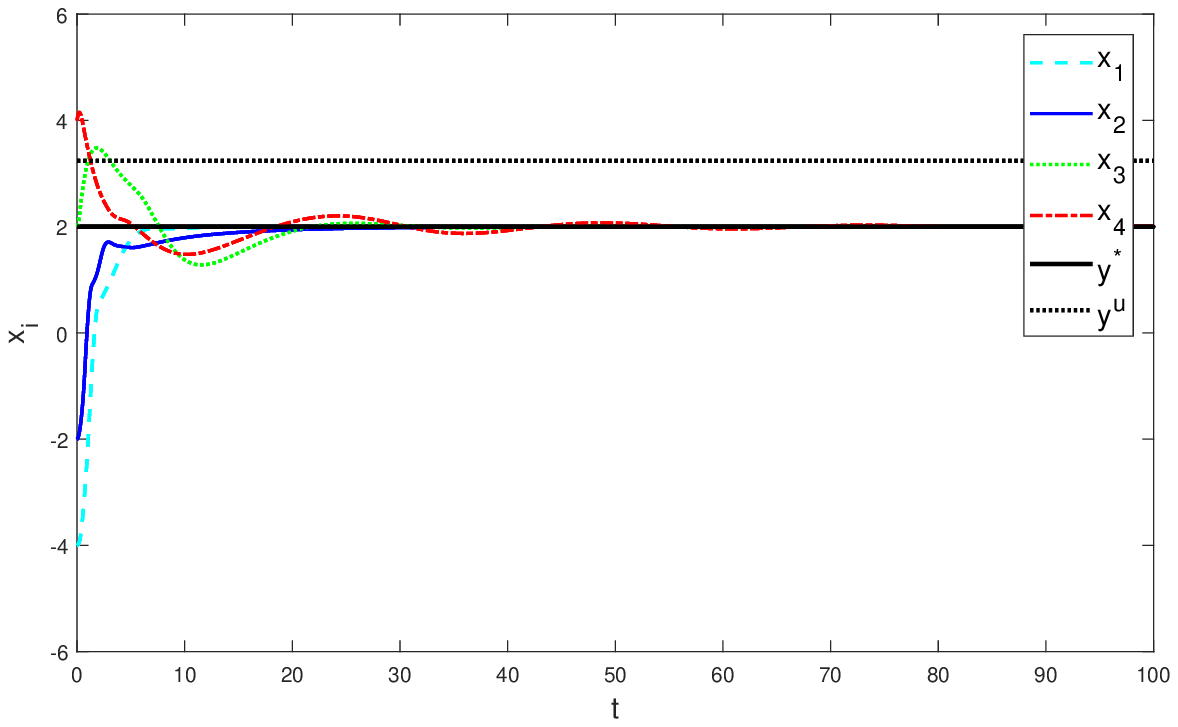}
\caption{Profiles of agent outputs in {\it   Example 2}.}	\label{fig:simu2-1} 
\end{figure}

{\it  Example 2}.  Consider another multi-agent system including two controlled Van der Pol oscillators 
\begin{align*}
\dot{x}_{i1}&=x_{i2}\\
\dot{x}_{i2}&=-(1+\mu_{i1})x_{i1}+(1+\mu_{i2})(\mu_{i3}-x_{i1}^2)x_{i2}+u_i+d_i, ~~ i=1,\,2
\end{align*}
and two controlled Duffing equations 
\begin{align*}
\dot{x}_{i1}&=x_{i2}\\
\dot{x}_{i2}&=-(1+\mu_{i1})x_{i1}(\mu_{i2}-x_{i1}^2)-(1+\mu_{i3})x_{i2}+u_i+d_i, ~~ i=3,\,4
\end{align*}
with input $u_i$, output $x_{i1}$, and disturbance $d_i$. Assume the disturbances are generated by \eqref{sys:disturbance} with 
\begin{align*}
D_i=1+\mu_{i4},\quad S_i=\begin{bmatrix} 0&1\\ -i&0
\end{bmatrix},\quad i=1,\,\dots,\,4
\end{align*}
and the unknown parameter $\mu_{ij}$ is randomly chosen between $-0.5$ and $0.5$. The information sharing graph is taken as the same with {\it  Example 1}.

\begin{figure}
\centering
\includegraphics[width=0.84\textwidth]{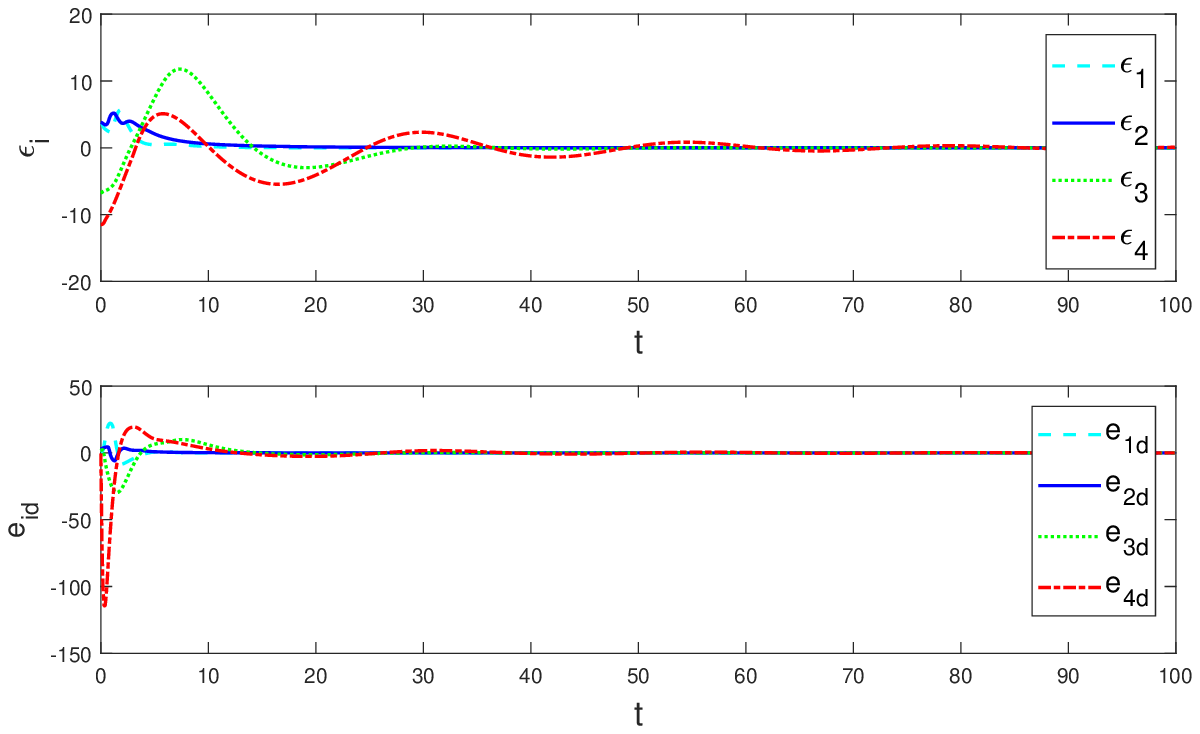}
\caption{Profiles of ${\bm{\epsilon}}_i$ and $e_{id}=u_{id}(t)+d_i(t)$  in {\it   Example 2}.}\label{fig:simu2-23}
\end{figure}

Although all agents are of the form \eqref{sys:agent} with $x_i=x_{i1}$, $n_i=2$, and $q=1$, the two classes of agent dynamics possess very different behaviors. This heterogeneity definitely brings many challenges in resolving their coordination problem. Moreover, we choose some complex cost functions as ${f_1}(y) = (y-8)^2$, ${f_2}(y) = \frac{y^2}{80\ln {({y^2} + 2} )} + (y - 5)^2$, ${f_3}(y)=\frac{y^2}{{20\sqrt {y^2 + 1} }} + y^2$, ${f_4}(y) = \ln \left( {{e^{ - 0.05{y}}} + {e^{0.05{y}}}} \right) + y^2$. 
with local interval constraints $[\,-3+i,\,1+i\,]$. 
Assumption \ref{ass:convex-strong} is fulfilled with $\underline{l}_i=1$ and $\bar{l}_i=3$ for $i=1,\,\dots,\, 4$.  Then, the formulated coordination problem for these agents can be solved by a controller of the form \eqref{ctr:full}. In fact, the optimal solution to the global constrained optimization problem is $y^*=2$ while the unconstrained optimal point is 	$y^{u}=3.24$ by directly minimizing $\sum_{i=1}^4 f_i(y)$.

%


For simulations, we choose the following matrices
\begin{align*}
F_i=\begin{bmatrix}
	-2i&1\\-i&0
\end{bmatrix},\quad G_i=\begin{bmatrix}
	-2i\\0
\end{bmatrix},\quad i=1,\,\dots,\,4
\end{align*} 
for the internal model \eqref{ctr:id} and use the same RBF neural network as in {\it Example 1}. The control gain function is chosen as $\rho_i(s)=s^6+1$ with parameters $k_{i1}=1$, $l=0.001$ for $1\leq i\leq 4$.  With randomly chosen initial conditions, the profiles of agent outputs under controller \eqref{ctr:full} are shown in Fig.~\ref{fig:simu2-1}.  We also list the approximation errors of feedforwarding input and external disturbance by neural networks and internal models in Fig.~\ref{fig:simu2-23}. It can be found that both $\e_i(t)$ and $e_{id}(t)$ converge towards zero as $t$ grows. 
The performance verifies the effectiveness of controller \eqref{ctr:full} to ensure the expected constrained optimal coordination for this heterogeneous uncertain multi-agent system \eqref{sys:agent}.

\section{Conclusion}\label{sec:con}
We have investigated the constrained optimal coordination problem for a class of heterogeneous nonlinear agents subject to both unknown dynamics and external disturbances. Jointly with internal model-based designs, novel distributed neural network-based controllers have been developed to overcome the technical difficulties brought by uncertainties, disturbances, and decision constraints under some standard assumptions.   Output feedback control for more general multi-agent systems and communication graphs will be our future work.   




\bibliographystyle{ieeetr}
\bibliography{opt-nonlinear-nn-abbr}

\end{document}